\documentclass[11pt,a4paper]{article}
\usepackage[T1]{fontenc}\usepackage{lmodern}
\ifdefined\XeTeXversion
\DeclareFontFamily{OT1}{lmr}{}\DeclareFontShape{OT1}{lmr}{m}{n}{<->rm-lmr10}{}\DeclareFontShape{OT1}{lmr}{bx}{n}{<->rm-lmbx10}{}
\DeclareFontFamily{OML}{lmm}{}\DeclareFontShape{OML}{lmm}{m}{it}{<-6>lmmi5<6-7>lmmi6<7-8>lmmi7<8-9>lmmi8<9->lmmi10}{}
\DeclareFontShape{OML}{lmm}{b}{it}{<-9>lmmib7<9->lmmib10}{}
\DeclareFontFamily{OMS}{lmsy}{}\DeclareFontShape{OMS}{lmsy}{m}{n}{<-7>lmsy6<7-8>lmsy7<8-9>lmsy8<9->lmsy10}{}
\DeclareFontFamily{OMX}{lmex}{}\DeclareFontShape{OMX}{lmex}{m}{n}{<->lmex10}{}
\DeclareFontFamily{T1}{lmr}{}
\DeclareFontShape{T1}{lmr}{m}{n}{<-7.5>ec-lmr6<7.5-9>ec-lmr8<9-10>ec-lmr9<10-12>ec-lmr10<12-16>ec-lmr12<16->ec-lmr17}{}
\DeclareFontShape{T1}{lmr}{m}{it}{<-10>ec-lmri9<10->ec-lmri10}{}
\DeclareFontShape{T1}{lmr}{bx}{n}{<-12>ec-lmbx10<12->ec-lmbx12}{}
\DeclareFontShape{T1}{lmr}{m}{sc}{<->ec-lmcsc10}{}
\DeclareFontFamily{T1}{lmtt}{}\DeclareFontShape{T1}{lmtt}{m}{n}{<->ec-lmtt10}{}
\fi
\usepackage{mathtools,amssymb,amsthm}
\usepackage{tikz}
\usepackage[margin=2.9cm]{geometry}
\usepackage{microtype}
\definecolor{linkblue}{rgb}{0,0,0.45}
\usepackage[bookmarksopen=true,bookmarksnumbered=true,colorlinks=true,
 allcolors=linkblue,pdfencoding=unicode]{hyperref}
\hypersetup{pdftitle={Absolute continuity of Dunkl intertwining measures via reflection chains},
 pdfauthor={Ziyang Jiang, Naihuan Jing, Ning Liu, Xingzhou Wang},
 pdfsubject={Representing measures of the Dunkl intertwining operator via reflection chains}}
\newtheorem{theorem}{Theorem}[section]
\newtheorem{corollary}[theorem]{Corollary}
\newtheorem{proposition}[theorem]{Proposition}
\newtheorem{lemma}[theorem]{Lemma}
\numberwithin{equation}{section}
\newcommand{\aaE}{\mathfrak a}
\newcommand{\R}{\mathbb R}
\newcommand{\C}{\mathbb C}
\newcommand{\E}{\mathbb E}
\newcommand{\spanR}{\operatorname{span}_{\R}}
\DeclareMathOperator{\supp}{supp}
\DeclareMathOperator{\ri}{ri}
\DeclareMathOperator{\conv}{conv}
\DeclareMathOperator{\aff}{aff}
\DeclareMathOperator{\Law}{Law}
\title{Absolute continuity of Dunkl intertwining measures
  \mbox{via reflection chains}}
\author{Ziyang Jiang\textsuperscript{1}, Naihuan Jing\textsuperscript{1},
  Ning Liu\textsuperscript{2}, and Xingzhou Wang\textsuperscript{3}\\[0.5em]
  \small\textsuperscript{1}Department of Mathematics, North Carolina State University, Raleigh, NC 27695, USA\\
  \footnotesize\href{mailto:zjiang29@ncsu.edu}{\nolinkurl{zjiang29@ncsu.edu}} (Z.~Jiang),
  \href{mailto:jing@ncsu.edu}{\nolinkurl{jing@ncsu.edu}} (N.~Jing)\\[0.25em]
  \small\textsuperscript{2}School of Mathematics, Hefei University of Technology, Hefei, Anhui 230009, China\\
  \footnotesize\href{mailto:mathliu123@outlook.com}{\nolinkurl{mathliu123@outlook.com}} (N.~Liu)\\[0.25em]
  \small\textsuperscript{3}Department of Statistics, North Carolina State University, Raleigh, NC 27695, USA\\
  \footnotesize\href{mailto:xwang268@ncsu.edu}{\nolinkurl{xwang268@ncsu.edu}} (X.~Wang)}
\date{}
\begin{document}
\maketitle
\begin{abstract}
Let $W$ be a finite real reflection group and $k\ge0$ a multiplicity
function. We prove that the representing measure $\mu_x$ of the Dunkl
intertwining operator is equivalent to Lebesgue measure on the relative
interior of its support. The support is the convex hull of the orbit of $x$
under the reflections with $k>0$. This proves a conjecture of R\"osler and
de Jeu. We also prove that the dual intertwining measures are absolutely
continuous when the roots with $k>0$ span the space.
\end{abstract}

\medskip
\noindent\textbf{Keywords.} Dunkl operators, intertwining operator,
representing measures, absolute continuity, reflection groups.

\noindent\textbf{2020 Mathematics Subject Classification.}
Primary 33C52. Secondary 33C67, 20F55, 42B10, 60J27.

\section{Introduction}\label{sec:intro}

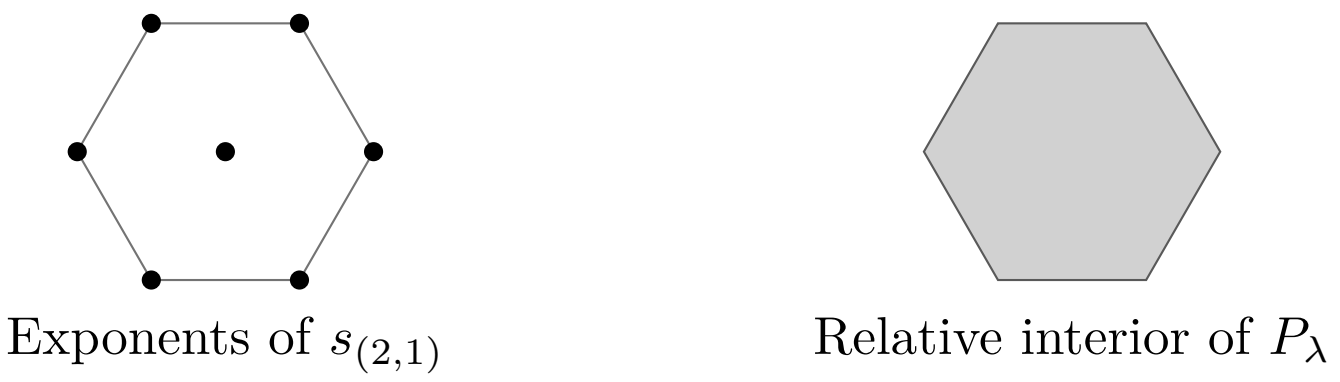
\begin{figure}[htbp]
\centering
\begin{tikzpicture}[x=1cm,y=1cm,font=\small]
  \def\hexpath{(0:0.98) -- (60:0.98) -- (120:0.98) --
               (180:0.98) -- (240:0.98) -- (300:0.98) -- cycle}
  \begin{scope}[xshift=-2.8cm]
    \draw[black!55] \hexpath;
    \foreach \angle in {0,60,120,180,240,300}
      \fill (\angle:0.98) circle (1.8pt);
    \fill (0,0) circle (1.8pt);
    \node[below] at (0,-0.96) {Exponents of $s_{(2,1)}$};
  \end{scope}
  \begin{scope}[xshift=2.8cm]
    \filldraw[fill=black!18,draw=black!65] \hexpath;
    \node[below] at (0,-0.96) {Relative interior of $P_\lambda$};
  \end{scope}
\end{tikzpicture}
\caption{The permutahedron for $\lambda=(2,1,0)$. The dots give the seven
exponent vectors of $s_{(2,1)}$. The shaded interior lies in the plane
$u_1+u_2+u_3=3$, which carries the area measure in the question above.}
\label{fig:hexagon}
\end{figure}

The transpositions in $S_n$ act on $\R^n$ as reflections. For $n=3$ the orbit of
$\lambda=(2,1,0)$ consists of the six vertices of the hexagon in
Figure~\ref{fig:hexagon}. The same hexagon is the Newton polytope of
the Schur polynomial $s_{(2,1)}=m_{(2,1)}+2m_{(1,1,1)}$ in three variables;
here $m_\nu$ denotes a monomial symmetric polynomial, so the central weight
occurs twice. For a partition $\lambda$ with at most $n$ parts, the Newton
polytope of $s_\lambda$ in $n$ variables is the permutahedron
$P_\lambda=\conv(S_n\lambda)$ \cite[Proposition~2.5]{r26}.

For constant multiplicity $k>0$, the transpositions $s_{ij}$ also give
the type-$A$ Dunkl operators
\[
 T_i f(u)=\partial_i f(u)+k\sum_{j\ne i}
             \frac{f(u)-f(s_{ij}u)}{u_i-u_j}.
\]
The intertwiner $V_k$, determined on graded polynomials by $V_k1=1$ and
$T_iV_k=V_k\partial_i$, gives a probability measure $\mu_\lambda^k$
on $P_\lambda$ representing $p\mapsto V_kp(\lambda)$ \cite{r17}.
Does $\mu_{(2,1,0)}^k$ have the same null sets as area on the hexagon?
For a finite real reflection group, we determine this measure class at
every $x$, including points on reflecting hyperplanes and multiplicities
that vanish on some reflections. This proves the conjecture of R\"osler
and de Jeu \cite[Remark~5]{r19}, formulated as Conjecture~A in \cite{r3}.

The measure is the law of an exponentially discounted occupation average of
a reflection chain on the orbit, and the finite-state Feynman--Kac formula
derives this positive integral representation directly from the kernel
equations. This gives a short probabilistic proof of R\"osler's positivity
theorem \cite[Theorems~1.1 and~1.2]{r17} for every finite real reflection
group, including noncrystallographic groups; the proofs of R\"osler and Voit
\cite{r20} and of McSwiggen and Sahi \cite[Remark~3.4]{r14} pass through the
trigonometric theory of crystallographic root systems. The chain also gives
the exact support $\conv(W_kx)$, where $W_k$ is generated by the reflections
of positive multiplicity; the general inclusion
$\supp\mu_x^k\subseteq\conv(Wx)$ is \cite[Corollary~5.3(i)]{r17}. Expanding
the kernel equation along rays, degree by degree, recovers the reflection-word
formula of Deleaval, Demni and Youssfi \cite{r5}. The holding times of the
chain provide barycentric coordinates on the orbit polytope, and these
determine the measure class; this step can also be derived from results of
Bouguet and Cloez \cite{r2} (see Section~\ref{sec:simplex}).

\subsection{Background}\label{sec:background}
Dunkl constructed the differential-reflection operators, their intertwiner
and its exponential kernel \cite{r7,r8}. Homogeneous expansions express the
moments of $V_k$ through products of reflections \cite{r8,r15,r16,r5}; the
reflection chain interprets these group-algebra formulas
(Section~\ref{sec:resolvents}), and its occupation law belongs to the class
of Markovian stick-breaking measures studied in \cite{r6}.

In type $A$ with $k>0$, Sawyer determined the support of each $\mu_x^k$ and
obtained a density on the fixed-sum hyperplane for the orbit average whenever
$x$ is not constant, including singular $x$
\cite[Theorem~4.2 and Corollary~4.2]{r21}. Analogous results hold for
$x\ne0$ in the $BC_q$ families with $d\in\{1,2,4\}$ and $p>2q-1$
\cite[Theorem~25 and Corollary~26]{r23}. The support argument uses the
$W$-invariance of $\supp\mu_x^k$ (see \cite{r9}). Since each individual
measure is dominated by $|W|$ times the orbit average, a density for the average gives absolute
continuity of each $\mu_x^k$. Theorem~\ref{thm:main} extends the exact
support description to arbitrary finite real reflection groups and gives
equivalence to relative Lebesgue measure for each individual measure. An
explicit density at regular points of type $A_2$ is given in
\cite[Corollary~2.2]{r1}, a rank-reduction formula in type $A$ in
\cite[Theorem~4.2]{r22}, and further formulas for dihedral groups in
\cite{r5,r4}; Maslouhi obtained densities for the Gaussian-smoothed operator
$V_k\circ e^{\Delta/2}$ \cite{r13}.

For a general reflection group, R\"osler and de Jeu proved absence of atoms
at regular $x$ (off all reflecting hyperplanes) when the total multiplicity
is positive, and asked for a density \cite[Theorem~3 and Remark~5]{r19}.
Brennecken and R\"osler formulated this as Conjecture~A, namely
$\mu_x^k\ll dx$ at every regular $x$ when the roots of positive multiplicity
span the ambient space, and identified gaps in the proof proposed by
Trim\`eche \cite{r25} for both the intertwiner and its dual
\cite[Section~1 and Section~2, items (3)--(5)]{r3}. Langen proved $L^2$
regularity at regular $x$ for spanning systems when all multiplicities exceed
$1/2$ \cite[Theorem~7.2]{r11}. The rank-one density shows that this uniform
$L^2$ threshold is sharp (Section~\ref{sec:rank-one}).

\subsection{Main results}\label{sec:main}

Let $\aaE$ be a finite-dimensional real Euclidean space and
$R\subset\aaE\setminus\{0\}$ a finite reduced root system, not necessarily
spanning $\aaE$ or crystallographic. Set
\[
 s_\alpha u=u-2\frac{\langle\alpha,u\rangle}{|\alpha|^2}\alpha,
 \qquad W=\langle s_\alpha:\alpha\in R\rangle .
\]
Choose a positive subsystem $R_+$ and a $W$-invariant multiplicity
$k:R\to[0,\infty)$. Write $k_\alpha=k(\alpha)$. The rational Dunkl
operators \cite{r7} are
\begin{equation}\label{eq:dunkl}
 T_\xi f(u)=\partial_\xi f(u)+
 \sum_{\alpha\in R_+}k_\alpha\langle\alpha,\xi\rangle
        \frac{f(u)-f(s_\alpha u)}{\langle\alpha,u\rangle}.
\end{equation}
\begingroup\setlength{\emergencystretch}{2em}
For smooth functions the quotients extend across the reflecting hyperplanes.
Nonreduced systems reduce to this case: the unit vectors $\beta$ on the
positive root rays, together with their negatives, form a reduced root system
with the $W$-invariant multiplicity
$\widetilde k_\beta=\sum_{\alpha\in R_+\cap\R_{>0}\beta}k_\alpha$. Since
$s_{c\beta}=s_\beta$ and
$\langle c\beta,\xi\rangle/\langle c\beta,u\rangle=\langle\beta,\xi\rangle/\langle\beta,u\rangle$
for $c>0$, neither the operators \eqref{eq:dunkl} nor the rates in
\eqref{eq:gen} change, and the weight $w_k$ of Section~\ref{sec:dual} changes
only by the positive constant $\prod_{\alpha\in R_{k,+}}|\alpha|^{2k_\alpha}$.
The intertwiner $V_k$ is the degree-preserving isomorphism of the graded
polynomial space satisfying $V_k1=1$ and $T_\xi V_k=V_k\partial_\xi$.
The Dunkl kernel $E_k$ is entire in both variables and is characterized by
\begin{equation}\label{eq:kernel}
 T_\xi^u E_k(u,z)=\langle\xi,z\rangle E_k(u,z),
 \qquad E_k(0,z)=1 .
\end{equation}
\par\endgroup
We extend the pairing complex bilinearly. In the algebraic construction
\cite{r8}, the degree-$m$ term of $E_k(\,\cdot\,,z)$ is
$V_k(\langle\,\cdot\,,z\rangle^m/m!)$.

Only the roots of positive multiplicity contribute in \eqref{eq:dunkl}.
Accordingly set
\[
 \begin{gathered}
 R_k=\{\alpha\in R:k_\alpha>0\},\qquad
 R_{k,+}=R_k\cap R_+,\qquad W_k=\langle s_\alpha:\alpha\in R_k\rangle,\\
 O_x=W_kx,\qquad K_x=\conv(O_x),\qquad H_x=\aff(O_x).
 \end{gathered}
\]
Since $k$ is $W$-invariant, $R_k$ is again a reduced root system, with
reflection group $W_k$. For an affine space $H$, let $\lambda_H$ denote its
Euclidean Lebesgue measure, with unit mass when $H$ is a point. Write $\ri$
for relative interior and $\sim$ for mutual absolute continuity.

\begin{theorem}\label{thm:main}
For every $x\in\aaE$ there is a unique compactly supported probability measure
$\mu_x^k$ such that, for $z\in\aaE_\C$ and every polynomial $p$,
\begin{equation}\label{eq:rep}
 E_k(x,z)=\int_{\aaE}e^{\langle y,z\rangle}\,d\mu_x^k(y),
 \qquad (V_kp)(x)=\int_{\aaE}p(y)\,d\mu_x^k(y).
\end{equation}
Its measure class and its support satisfy
\begin{equation}\label{eq:class}
 \mu_x^k\sim\mathbf 1_{\ri K_x}\lambda_{H_x},
 \qquad \supp\mu_x^k=K_x .
\end{equation}
In particular,
\begin{equation}\label{eq:fulldim}
 \mu_x^k\ll\lambda_{\aaE}\quad\Longleftrightarrow\quad H_x=\aaE .
\end{equation}
The measure has the realization
\begin{equation}\label{eq:chain}
 \mu_x^k=\Law\!\left(\int_0^\infty e^{-s}X_s\,ds\right),
\end{equation}
where $X_0=x$ and the chain on $O_x$ has generator
\begin{equation}\label{eq:gen}
 (Qf)(v)=\sum_{\alpha\in R_{k,+}}k_\alpha\bigl(f(s_\alpha v)-f(v)\bigr).
\end{equation}
For $k=0$ this gives $\mu_x^k=\delta_x$.
\end{theorem}

For $\lambda=(2,1,0)$ and $k>0$, \eqref{eq:class} answers the opening
question: $\mu_\lambda^k$ has a density with respect to area on the hexagon
that is positive almost everywhere, and the six edges carry no mass.

\begin{corollary}\label{cor:positive}
For every real polynomial $p$ and every $x\in\aaE$,
\[
 p|_{K_x}\ge0\ \Longrightarrow\ V_kp(x)\ge0 .
\]
For $r>0$, put $B_r=\{x:|x|\le r\}$ and let $f=\sum_{m\ge0}f_m$,
where $f_m$ is a homogeneous polynomial of degree $m$. Then
\[
 \sum_{m\ge0}\sup_{B_r}|f_m|<\infty
 \quad\Longrightarrow\quad
 \sum_{m=0}^{N}V_kf_m(x)\xrightarrow[N\to\infty]{}
 V_kf(x):=\int_\aaE f\,d\mu_x^k
\]
uniformly for $x\in B_r$.
\end{corollary}

\begin{corollary}\label{cor:hull}
Let $R_1,\ldots,R_m$ be the orthogonal irreducible components of $R_k$ and
$E_j=\spanR R_j$. Take $m=0$ if $R_k$ is empty. Write
$\aaE=E_0\oplus E_1\oplus\cdots\oplus E_m$, where
$E_0=(\spanR R_k)^\perp$, and let $x_j$ be the projection of $x$ onto $E_j$.
Then
\begin{equation}\label{eq:hull}
 H_x=x_0+\bigoplus_{\substack{1\le j\le m\\x_j\ne0}}E_j .
\end{equation}
Consequently $\mu_x^k\ll\lambda_{\aaE}$ precisely when $E_0=\{0\}$ and
$x_j\ne0$ for every $1\le j\le m$.
\end{corollary}

In particular, Conjecture~A (Section~\ref{sec:background}) holds:
\[
 \spanR R_k=\aaE,\ \ \langle\alpha,x\rangle\ne0\ (\alpha\in R)
 \ \Longrightarrow\ E_0=\{0\},\ \ x_j\ne0\ (1\le j\le m)
 \ \Longrightarrow\ \mu_x^k\ll\lambda_\aaE .
\]
The first implication holds because $x_j=0$ would give
$\langle\alpha,x\rangle=\langle\alpha,x_j\rangle=0$ for $\alpha\in R_j$.
If $R_k$ is irreducible and spans $\aaE$, then $m=1$, $E_0=\{0\}$ and
$x_1=x$, so $\mu_x^k\ll\lambda_\aaE$ for every $x\ne0$, including points on
the walls.

\begin{theorem}\label{thm:dual}
Suppose that the roots of positive multiplicity span $\aaE$. Let
$(\nu_y^k)_{y\in\aaE}$ be the vaguely continuous family representing the
dual intertwining operator $\prescript{t}{}{V_k}$, as in
Lemma~\ref{lem:dual-smooth}. Then $\nu_y^k\ll\lambda_{\aaE}$ for every $y\in\aaE$.
\end{theorem}

Section~\ref{sec:dual} proves Theorem~\ref{thm:dual} in steps (a) and (b).
Steps (c) and (d) give a quantitative bound near the walls and treat the
case where $R_k$ does not span $\aaE$.
\begin{enumerate}
\item[(a)] \emph{Almost every $y$.} By the transpose formula
\eqref{eq:transpose}, Theorem~\ref{thm:main} gives
$\nu_y^k\ll\lambda_\aaE$ for almost every $y$.
\item[(b)] \emph{A given $y$.} We move a cutoff of a null set along a ray to
a point where absolute continuity is already known. The reflection terms near
the walls are nonnegative, because the cutoff decreases in the directions
normal to the walls.
\item[(c)] \emph{Mass near the walls.} For the local cutoffs $f_\varepsilon$
at a reflection stratum, Corollary~\ref{cor:wall} gives
$\sup_{y\in\aaE}\nu_y^k(f_\varepsilon)=O(\varepsilon^{2\gamma_*})$,
where $\gamma_*=\sum_{\alpha\in J} k_\alpha$ for the positive active roots
$J$ vanishing on the stratum. The proof uses the Fourier bound of
Proposition~\ref{prop:dual-wall}.
\item[(d)] \emph{The nonspanning case.} If $R_k$ does not span $\aaE$, then
by \eqref{eq:factor} the measure $\nu_y^k$ is absolutely continuous on the
affine subspace $y+\spanR R_k$.
\end{enumerate}

The paper is organized as follows. Section~\ref{sec:words} proves
\eqref{eq:rep}, \eqref{eq:chain} and Corollary~\ref{cor:positive} and relates
them to reflection words. Section~\ref{sec:simplex} proves \eqref{eq:class},
\eqref{eq:fulldim} and Corollary~\ref{cor:hull}. Section~\ref{sec:dual} treats
the dual measures, and Section~\ref{sec:examples} computes moments and
discusses the permutahedron, the Birkhoff polytope and the rank-one case.

\section{The reflection chain}\label{sec:words}

Fix \(x\in\aaE\), \(z\in\aaE_\C\) and \(O=W_kx\). On functions on \(O\) put
\(d_z(v)=\langle v,z\rangle\), \(D_zf=d_zf\),
\(Cf(v)=\sum_{\alpha\in R_{k,+}}k_\alpha f(s_\alpha v)\) and
\(\gamma=\sum_{\alpha\in R_{k,+}}k_\alpha\), so \(Q=C-\gamma I\).
Start \(G_s\) at \(e\in W_k\), with jumps \(g\mapsto s_\alpha g\) at rate
\(k_\alpha\). Then \(X_s=G_sx\) has generator \(Q\) and the same random
linear map \(A\) serves all starting points:
\begin{equation}\label{eq:A}
 A=\int_0^\infty e^{-s}G_s\,ds,\qquad \|A\|\le1,\qquad
 Ax=\sum_{v\in O}\left(\int_0^\infty e^{-s}\mathbf1_{\{X_s=v\}}\,ds\right)v\in K_x .
\end{equation}
Indeed each \(G_s\) is orthogonal and the occupation weights sum to one.
The Dunkl-theoretic input is \eqref{eq:kernel}, \(T_\xi V_k=V_k\partial_\xi\),
\(V_k1=1\) and Dunkl's construction~\cite{r8}.

\subsection{The Feynman--Kac formula}
\begin{proof}[Proof of \eqref{eq:rep}, \eqref{eq:chain} and Corollary~\ref{cor:positive}]
For \(F_v(t)=E_k(tv,z)\), the kernel equation and \eqref{eq:dunkl} give
\begin{equation}\label{eq:ray}
 \begin{aligned}
 t\,d_z(v)F_v(t)&=tT_v^uE_k(tv,z)=tF_v'(t)+\sum_{\alpha\in R_{k,+}}k_\alpha
       \bigl(F_v(t)-F_{s_\alpha v}(t)\bigr)
 =tF_v'(t)-(QF(t))(v).
 \end{aligned}
\end{equation}
For \(t>0\), \(t\langle\alpha,v\rangle=\langle\alpha,tv\rangle\) cancels the
denominator whenever \(s_\alpha v\ne v\); the other summands vanish.
Thus \(tF'=QF+tD_zF\), \(F(0)=\mathbf1\). With \(h(s,v)=F_v(e^{-s})\),
\(I_s=\int_0^s e^{-u}d_z(X_u)\,du\) and \(M_s=e^{I_s}h(s,X_s)\), we obtain
\begin{equation}\label{eq:FK}
 \begin{gathered}
 (\partial_s+Q)h=-e^{-s}D_zh,\qquad
 \operatorname{drift}(M_s)=e^{I_s}\bigl[(\partial_s+Q)h+e^{-s}D_zh\bigr](s,X_s)=0,\\
 E_k(x,z)=M_0=\lim_{s\to\infty}\E M_s=\E e^{\langle Ax,z\rangle}.
 \end{gathered}
\end{equation}
The drift formula is the finite-state chain rule. The bound
\(|M_s|\le e^{|x||z|}\max_{v\in O,\,0\le t\le1}|F_v(t)|\) makes \(M\) a
martingale and permits dominated convergence, since \(h(s,\cdot)\to\mathbf1\)
and \(I_s\to\langle Ax,z\rangle\). Define \(\mu_x^k=\Law(Ax)\). Taking \(z=i\xi\)
in \eqref{eq:FK} gives uniqueness by the uniqueness theorem for characteristic
functions. For \(p_n(u)=\langle u,z\rangle^n/n!\),
\[
 (V_kp_n)(x)=[t^n]E_k(tx,z)=[t^n]\E e^{t\langle Ax,z\rangle}
 =\E p_n(Ax)=\int_\aaE p_n\,d\mu_x^k .
\]
Here \(E_k(\cdot,z)=\sum_{m\ge0}V_kp_m\) is the homogeneous expansion,
\(A(tx)=tAx\), and \(|\langle Ax,z\rangle|\le |x||z|\).
Polarization gives \(V_kp(x)=\E p(Ax)\) for every polynomial \(p\).
By \eqref{eq:A}, \(p|_{K_x}\ge0\) gives \(V_kp(x)\ge0\), and
\(\sup_{|x|\le r}|V_kp(x)|\le\sup_{|y|\le r}|p(y)|\). Hence
\(\sum_m\sup_{|x|\le r}|V_kf_m(x)|
\le\sum_m\sup_{|y|\le r}|f_m(y)|<\infty\), so uniform passage to the limit
in the integral formula gives Corollary~\ref{cor:positive}.
The same bound and uniform limit hold for bounded Borel summands on \(B_r\)
when \(V_kg(x)=\int_{B_r} g\,d\mu_x^k\) for \(x\in B_r\).
\end{proof}

\subsection{Resolvents and reflection words}\label{sec:resolvents}

Put \(h_m=[t^m]F=(V_kp_m)|_O\). Comparing coefficients in \eqref{eq:ray}
gives \(h_0=\mathbf1\), \((mI-Q)h_m=D_zh_{m-1}\) and, for \(m\ge1\),
\begin{equation}\label{eq:rec}
 \begin{aligned}
 R_m&=(mI-Q)^{-1}=\int_0^\infty e^{-ms}e^{sQ}\,ds
       =\sum_{\ell\ge0}\frac{C^\ell}{(m+\gamma)^{\ell+1}},\\
 h_m&=R_mD_zh_{m-1}=(R_mD_z)\cdots(R_1D_z)\mathbf1 .
 \end{aligned}
\end{equation}
Herbach derives the corresponding ray equation and coefficient recursion
from the stationary joint Laplace transform
\cite[Proposition~4.2 and Corollary~4.3]{r10}.
The Markov semigroup gives the Laplace integral and
\(\|R_m\|_{\infty\to\infty}\le1/m\); \(\|C\|_{\infty\to\infty}\le\gamma\)
gives the Neumann series, and \(\|D_z\|_{\infty\to\infty}\le|x||z|\)
gives \(\|h_m\|_\infty\le(|x||z|)^m/m!\). Compare the ray system in
\cite[Section~3]{r19} and the Euler inverse in \cite[Section~2]{r8}.

On \(W_k\) let \(C_W\varphi(g)=\sum_{\alpha\in R_{k,+}}k_\alpha\varphi(gs_\alpha)\)
and \(Q_W=C_W-\gamma I\). Conjugation preserves the weighted reflections,
so \(Q_W(f(\cdot\,x))=(Qf)(\cdot\,x)\). For the total weight \(c_\ell(g)\)
of reflection words of length \(\ell\) with product \(g\), we have
\[
 \begin{aligned}
 c_\ell(g)&=\sum_{s_{\alpha_1}\cdots s_{\alpha_\ell}=g}
                 k_{\alpha_1}\cdots k_{\alpha_\ell},\qquad
 c_0=\mathbf1_{\{e\}},\qquad (C_W^\ell)_{g,g'}=c_\ell(g^{-1}g'),\\
 [(mI-Q_W)^{-1}]_{g,g'}&=\sum_{\ell\ge0}
     \frac{c_\ell(g^{-1}g')}{(m+\gamma)^{\ell+1}}=:C_m(g^{-1}g').
 \end{aligned}
\]
The word sum is over \(\alpha_i\in R_{k,+}\); the second line is
\eqref{eq:rec} on \(W_k\). With \(g_0=e\), the two moment formulas now read
\begin{equation}\label{eq:words}
 \begin{aligned}
 \frac{V_k(\langle\cdot,z\rangle^n)(x)}{n!}
 &=\frac{\E\langle Ax,z\rangle^n}{n!}
 =\int_{0<t_1<\cdots<t_n}e^{-\sum_jt_j}\,
           \E_x\!\prod_{j=1}^n d_z(X_{t_j})\,dt\\
 &=\bigl[(R_nD_z)\cdots(R_1D_z)\mathbf1\bigr](x)\\
 &=\sum_{g_1,\ldots,g_n\in W_k}\prod_{j=1}^n
       \bigl(C_{n+1-j}(g_{j-1}^{-1}g_j)\langle g_jx,z\rangle\bigr).
 \end{aligned}
\end{equation}
For a resolvent formula for mixed moments of the simplex occupation
weights, see \cite[Proposition~3.8]{r12}.
The first line uses \eqref{eq:rep} and orders the \(n\) observation times
in \eqref{eq:A}. For \(t_0=0\), the Markov property and
\(\sum_j t_j=\sum_j(n+1-j)(t_j-t_{j-1})\) assign \(R_{n+1-j}\) to the
waiting time before observation \(j\). The last line multiplies the lifted
kernels. Reversing the indices \(g_j\mapsto g_{n+1-j}\) gives
\cite[Proposition~3.1 and Corollary~3.2]{r5}; inactive roots have weight
zero and the convergence condition is \(\sum_{\alpha\in R_+}|k_\alpha|
=\gamma<1+\gamma=|1+\gamma|\).

\section{Measure class and affine hull}\label{sec:simplex}

\begin{proposition}\label{prop:simplex}
Let \(O\) be a nonempty finite subset of a Euclidean space,
\(K=\conv O\), and \(H=\aff O\). For an irreducible time-homogeneous continuous-time Markov chain on \(O\) with \(X_0=v\in O\),
\[
 Y=\int_0^\infty e^{-s}X_s\,ds,\qquad
 \Law_v(Y)\sim\mathbf1_{\ri K}\lambda_H,\qquad \supp\Law_v(Y)=K .
\]
\end{proposition}

\begin{proof}
If \(|O|=1\), then \(Y=v\). Otherwise let \(V_j,T_j,H_j=T_{j+1}-T_j\)
be the successive states, jump times (\(T_0=0\)) and holding times.
Irreducibility gives \(N=\min\{n:\{V_0,\ldots,V_{n-1}\}=O\}<\infty\)
almost surely. Given \(\mathcal G=\sigma((V_j)_{j\ge0},(H_{N+j})_{j\ge0})\),
\(N\) and \(z=\int_0^\infty e^{-s}X_{T_N+s}\,ds\in K\) are fixed, and
\(H_0,\ldots,H_{N-1}\) are independent exponentials with positive rates.
Thus \(u=(e^{-T_1},\ldots,e^{-T_N})\) has a positive conditional density on
\(D=\{1>u_1>\cdots>u_N>0\}\): with \(u_0=1\), the inverse change of
variables is \(H_j=\log(u_j/u_{j+1})\). Splitting at the jumps gives
\[
 \begin{aligned}
 Y&=\sum_{j=0}^{N-1}(u_j-u_{j+1})V_j+u_Nz=:L(u),\\
 L(D)&=\ri L(\overline D)=\ri\conv\{V_0,\ldots,V_{N-1},z\}=\ri K,\\
 \Law_v(Y\mid\mathcal G)&=L_*\Law_v(u\mid\mathcal G)
       \sim L_*(\mathbf1_D\,du)\sim\mathbf1_{\ri K}\lambda_H .
 \end{aligned}
\]
The coefficients range over the closed barycentric simplex on \(\overline D\),
and affine maps preserve relative interiors \cite[Theorem~6.6]{rockafellar}.
Fubini gives the last equivalence: each fiber over \(\ri K\) meets \(D\)
in a nonempty relatively open set of positive fiber measure.
Averaging over \(\mathcal G\) yields the asserted equivalence and
\(\supp\Law_v(Y)=\overline{\ri K}=K\).
\end{proof}

\begin{proof}[Completion of the proof of Theorem~\ref{thm:main}]
The positive-rate reflections generate \(W_k\), so the chain on \(O_x\) is
irreducible. Proposition~\ref{prop:simplex} and \eqref{eq:chain} give
\eqref{eq:class}, hence \(\mu_x^k\ll\lambda_\aaE\) if \(H_x=\aaE\);
otherwise \(\lambda_\aaE(H_x)=0<1=\mu_x^k(H_x)\). This proves \eqref{eq:fulldim}.
\end{proof}

Alternatively, Proposition~\ref{prop:simplex} follows from \cite{r2}. Take a
two-sided stationary chain $I$ with the time-reversed generator, let \(e_w\) denote the vertices of the probability simplex \(\Delta_O\),
and set \(\pi(p)=\sum_{w\in O}p_ww\). The stationary solution of
\(\dot Z_s=e_{I_s}-Z_s\) gives
\[
 Z_0=\int_{-\infty}^0 e^s e_{I_s}\,ds,\qquad
 \Law_v(Y)=\pi_*\Law(Z_0\mid I_0=v)\sim\mathbf1_{\ri K}\lambda_H .
\]
Indeed \(I_{-s}\), conditional on \(I_0=v\), has the original generator.
The conditional law on \(\Delta_O\) is equivalent to relative Lebesgue
measure on its interior by the minorization and absolute continuity in
\cite[Lemma~3.2 and Proposition~3.4]{r2}; the same affine-fiber argument
applies to \(\pi\). For a stick-breaking representation of the joint stationary law of
\((I_s,Z_s)\), see \cite[Theorem~3.7]{r12}.

\begin{proof}[Proof of Corollary~\ref{cor:hull}]
Put \(W_j=\langle s_\alpha:\alpha\in R_j\rangle\) and \(O_j=W_jx_j\).
Since \(W_k=\prod_{j=1}^m W_j\) acts componentwise and fixes \(E_0\),
\[
H_x=\aff(W_kx)=x_0+\aff\Bigl(\prod_{j=1}^m O_j\Bigr)=x_0+\bigoplus_{j=1}^m\aff O_j=x_0+\bigoplus_{j=1}^m\spanR O_j=x_0+\bigoplus_{\substack{1\le j\le m\\ x_j\ne0}}E_j .
\]
Here affine hulls commute with products, and
\(\frac1{|W_j|}\sum_{w\in W_j}wx_j\in E_j^{W_j}=\{0\}\)
gives \(0\in\aff O_j\), hence \(\aff O_j=\spanR O_j\).
By irreducibility of \(E_j\), the \(W_j\)-invariant span \(\spanR O_j\)
equals \(E_j\) when \(x_j\ne0\) and \(\{0\}\) when \(x_j=0\).
Thus \eqref{eq:fulldim} gives
\[
 \mu_x^k\ll\lambda_{\aaE}
 \quad\Longleftrightarrow\quad H_x=\aaE
 \quad\Longleftrightarrow\quad
 E_0=\{0\}\ \text{and}\ x_j\ne0\ \text{for every }1\le j\le m .
 \qedhere
\]
\end{proof}

Assume all multiplicities are positive. For $W=S_3$, put
$x=(1,1,-2)$ and $H=\{u\in\R^3:u_1+u_2+u_3=0\}$. Then
\[
 \begin{aligned}
 (12)x&=x, \qquad
 Wx=\{(1,1,-2),(1,-2,1),(-2,1,1)\},\\
 H_x&=H,\qquad \dim K_x=2,\qquad
 \mu_x^k\sim\mathbf1_{\ri K_x}\lambda_H .
 \end{aligned}
\]
For $W=A_1\times A_1$ acting by coordinate sign changes and
$y=(a,0)$ with $a>0$,
\[
 Wy=\{(-a,0),(a,0)\},\qquad K_y=[-a,a]\times\{0\},\qquad
 \mu_y^k\sim\bigl(\mathbf1_{(-a,a)}(t)\,dt\bigr)\otimes\delta_0 .
\]

\section{Absolute continuity of the dual measures}\label{sec:dual}

Let $d=\dim\aaE$. By Theorem~\ref{thm:main}, $V_k$ extends to continuous
functions by $V_kg(x)=\int_\aaE g\,d\mu_x^k$, and with the weight
$w_k(x)=\prod_{\alpha\in R_{k,+}}|\langle\alpha,x\rangle|^{2k_\alpha}$ the dual
intertwiner is the weighted transpose
\begin{equation}\label{eq:transpose}
 \int_{\aaE}\prescript{t}{}{V_k}f(y)g(y)\,dy
   =\int_{\aaE} f(x)V_kg(x)w_k(x)\,dx
 \qquad (f,g\in C_c^\infty(\aaE)).
\end{equation}
Trim\`eche~\cite{r24} showed that $\prescript{t}{}{V_k}$ is represented by
measures. The following lemma derives the properties we need directly from
Theorem~\ref{thm:main}.

\begin{lemma}\label{lem:dual-smooth}
For $f\in C_c^\infty(\aaE)$ and $v\in\aaE$ put
\[
 \mathcal F_kf(\xi)=\int_\aaE f(x)E_k(x,-i\xi)w_k(x)\,dx,\qquad
 \nu^k_y(f)=(2\pi)^{-d}\int_\aaE e^{i\langle y,\xi\rangle}\mathcal F_kf(\xi)\,d\xi .
\]
Then $\mathcal F_k(T_vf)(\xi)=i\langle v,\xi\rangle\mathcal F_kf(\xi)$,
$\mathcal F_kf$ is rapidly decreasing, and $\nu^k_y(f)$ is the smooth
representative of $\prescript{t}{}{V_k}f$ in \eqref{eq:transpose}; it is the
integral of $f$ against a unique positive Radon measure $\nu^k_y$, and
\begin{equation}\label{eq:diffdual}
 \partial_{y,v}\nu_y^k(f)=\nu_y^k(T_vf).
\end{equation}
The family is vaguely continuous, and $\sup_{y\in Y}\nu^k_y(K)<\infty$ for
compact $Y,K\subset\aaE$.
\end{lemma}

\begin{proof}
Write $\beta_f(dy)=\int_\aaE f(x)w_k(x)\mu_x^k(dy)\,dx$.
Weighted integration by parts \cite[Proposition~2.1]{r18},
\eqref{eq:kernel} and \eqref{eq:FK} give, for $N\ge0$ and $\xi\ne0$,
\[
 \begin{aligned}
 \mathcal F_k(T_v f)(\xi)
 &=-\int_\aaE f(x)\,T_v^x E_k(x,-i\xi)w_k(x)\,dx
   =i\langle v,\xi\rangle\,\mathcal F_k f(\xi),\\
 |\xi|^N|\mathcal F_kf(\xi)|
 &=\bigl|\mathcal F_k(T_{\xi/|\xi|}^Nf)(\xi)\bigr|
   \le\sup_{|u|=1}\int_\aaE|T_u^Nf|\,w_k\,dx<\infty .
 \end{aligned}
\]
In the first line one may cut off the kernel invariantly near $W_k\supp f$.
The second uses $|E_k(x,-i\xi)|\le1$ and smooth divided differences:
$T_u^N f$ has support in $W_k\supp f$ and depends polynomially on $u$.
Since $\widehat{\beta_f}=\mathcal F_k f$ by \eqref{eq:rep}, Fourier
inversion, Fubini and the multiplier identity give
\[
 \begin{aligned}
 \int_\aaE g(y)\nu_y^k(f)\,dy
 &=\int_\aaE g\,d\beta_f=\int_\aaE f\,V_kg\,w_k\,dx,\\
 \partial_{y,v}\nu_y^k(f)
 &=(2\pi)^{-d}\int_\aaE e^{i\langle y,\xi\rangle}
        i\langle v,\xi\rangle\,\mathcal F_k f(\xi)\,d\xi=\nu_y^k(T_v f).
 \end{aligned}
\]
Rapid decay permits differentiation to every order. For $f\ge0$ the
smooth density $\nu_\cdot^k(f)$ of $\beta_f\ge0$ is nonnegative everywhere.
In particular $|\nu_y^k(f)|\le\|f\|_\infty\nu_y^k(\theta)$ for
$\theta\in C_c^\infty$ with $\theta\ge\mathbf1_{\supp f}$, so $\nu_y^k$
extends to a unique positive Radon measure. For compact $K,Y$ take
$\theta\ge\mathbf1_K$.
Then $\sup_{y\in Y}\nu_y^k(K)\le\sup_{y\in Y}\nu_y^k(\theta)<\infty$,
and uniform approximation on a common compact support gives vague continuity.
\end{proof}

\begin{proof}[Proof of Theorem~\ref{thm:dual}]
The case \(d=0\) is immediate. Assume \(d>0\).

\smallskip
\noindent\emph{(a) Almost every \(y\).}
Since \(\spanR R_k=\aaE\), Corollary~\ref{cor:hull} gives
\(\mu_x^k\ll dy\) outside the finite union of proper subspaces
\(\bigcup_{j=1}^m\{x_j=0\}\). Put
\(M(dx,dy)=w_k(x)\,dx\,\mu_x^k(dy)\). For \(f,g\in C_c^\infty(\aaE)\),
\[
 \int g(y)\nu_y^k(f)\,dy
 \overset{\eqref{eq:transpose}}{=}\int f(x)V_kg(x)w_k(x)\,dx
 =\int f(x)g(y)\,M(dx,dy).
\]
Thus \(M(dx,dy)=\nu_y^k(dx)\,dy\). If \(Z\subset\aaE^2\) is Borel
and null, then its sections \(Z_x=\{y:(x,y)\in Z\}\) are null for
a.e. \(x\), so \(M(Z)=\int w_k(x)\mu_x^k(Z_x)\,dx=0\).
Writing \(M(dx,dy)=p(x,y)\,dx\,dy\), we obtain, for every bounded
half-open box \(U\) with rational endpoints and every bounded Borel \(B\),
\[
 \int_B\nu_y^k(U)\,dy
 =M(U\times B)
 =\int_B\!\left(\int_U p(x,y)\,dx\right)dy .
\]
Since these \(U\) form a countable determining family,
\(\nu_y^k(dx)=p(x,y)\,dx\ll dx\) for a.e. \(y\).

\smallskip
\noindent\emph{(b) A given \(y\).}
Fix \(y_0\in\aaE\). By inner regularity and the finite stratification
by \(J(x)=\{\alpha\in R_{k,+}:\langle\alpha,x\rangle=0\}\), it suffices
to show \(\nu_{y_0}^k(N)=0\) for every compact null set \(N\) on which
\(J(x)=J\) is constant. Let \(P\) be the orthogonal projection onto
\(L=\spanR J\) and put \(P^\perp=I-P\), so \(N\subset L^\perp\).
Choose \(\delta>0\) so that the closed \(3\delta\)-neighborhood \(K\)
of \(N\) misses \(\alpha^\perp\) for every
\(\alpha\in R_{k,+}\setminus J\). By (a) and polar coordinates about
\(y_0\), there exist a unit vector \(v\) and \(0<\ell<\delta\) such that,
with \(y_t=y_0+tv\),
\[
 \nu_{y_\ell}^k\ll dx,\qquad
 \nu_{y_t}^k\ll dx\quad\text{for a.e. }t\in(0,\ell).
\]

Choose \(0\le a_n\le1\) in \(C_c^\infty(L^\perp)\), equal to one on
\(N\) and supported within distance \(1/n\) of \(N\), and choose
\(0\le\psi\in C_c^\infty(\R)\) with \(\psi(0)=1\) and
\(\psi'\le0\) on \([0,\infty)\). Put
\[
 f_n(x)=a_n(P^\perp x)\psi(n^2|Px|^2),\qquad
 P_tf(x)=f(x-tv),\qquad
 T_v=\partial_v+\sum_{\alpha\in R_{k,+}}B_{\alpha,v},
\]
where \(B_{\alpha,v}\) is the \(\alpha\)-summand of \eqref{eq:dunkl}.
Then \(\mathbf1_N\le f_n\le1\), \(f_n\to\mathbf1_N\) pointwise, and
\(\supp P_tf_n\subset K\) for \(0\le t\le\ell\) and all sufficiently
large \(n\).

For \(\alpha\in J\), put \(r=n^2|P(x-tv)|^2\) and
\(r'=n^2|P(s_\alpha x-tv)|^2\). Since \(P^\perp s_\alpha x=P^\perp x\),
\[
 \begin{aligned}
 r'-r
 &=\frac{4tn^2\langle\alpha,x\rangle\langle\alpha,v\rangle}{|\alpha|^2},\\
 B_{\alpha,v}P_tf_n(x)
 &=k_\alpha a_n(P^\perp(x-tv))
   \frac{\langle\alpha,v\rangle}{\langle\alpha,x\rangle}
   [\psi(r)-\psi(r')]\\
 &=-\frac{4tn^2k_\alpha\langle\alpha,v\rangle^2}{|\alpha|^2}
   a_n(P^\perp(x-tv))
   \int_0^1\psi'((1-s)r+sr')\,ds\ge0 .
 \end{aligned}
\]
The last expression extends continuously to \(\alpha^\perp\).
For \(\alpha\in R_{k,+}\setminus J\), the denominator
\(|\langle\alpha,x\rangle|\) is bounded away from zero on
\(K\cup s_\alpha K\), hence
\(\lvert B_{\alpha,v}P_tf_n\rvert\le C\mathbf1_{K\cup s_\alpha K}\).

For each fixed \(n\), the difference quotients of \(P_t f_n\)
with respect to \(t\) converge uniformly on a common compact support. Thus \eqref{eq:diffdual}, vague continuity and the local
mass bound give the product rule in the first line below. Integrating
and using the sign of the terms indexed by \(J\), we obtain
\[
 \begin{aligned}
 \frac d{dt}\nu_{y_t}^k(P_tf_n)
 &=\nu_{y_t}^k\bigl((T_v-\partial_v)P_tf_n\bigr)
  =\sum_{\alpha\in R_{k,+}}\nu_{y_t}^k(B_{\alpha,v}P_tf_n),\\
 0\le\nu_{y_0}^k(N)\le\nu_{y_0}^k(f_n)
 &=\nu_{y_\ell}^k(P_\ell f_n)
   -\int_0^\ell\sum_{\alpha\in R_{k,+}}
      \nu_{y_t}^k(B_{\alpha,v}P_tf_n)\,dt\\
 &\le\nu_{y_\ell}^k(P_\ell f_n)
   +\sum_{\alpha\in R_{k,+}\setminus J}
      \int_0^\ell\nu_{y_t}^k(|B_{\alpha,v}P_tf_n|)\,dt
   \longrightarrow0 .
 \end{aligned}
\]
Indeed, \(P_tf_n\to\mathbf1_{N+tv}\), and both \(N+tv\) and
\(s_\alpha(N+tv)\) are null. Absolute continuity at \(t=\ell\) and
for a.e. \(t\in(0,\ell)\) gives dominated convergence in \(x\),
and the uniform mass bound on \(K\cup s_\alpha K\) gives dominated
convergence in \(t\). Hence \(\nu_{y_0}^k\ll dx\) for every \(y_0\in\aaE\).
\end{proof}

\subsection{Mass near the walls}\label{sec:wall}
This subsection carries out step~(c) of Section~\ref{sec:main}.

\begin{proposition}\label{prop:dual-wall}
For $x_*$ on an active hyperplane put
\[
 J=\{\alpha\in R_{k,+}:\langle\alpha,x_*\rangle=0\},\quad
 L=\spanR J,\quad r=\dim L,\quad\gamma_*=\sum_{\alpha\in J}k_\alpha>0 ,
\]
let $P,P^\perp$ be the orthogonal projections onto $L,L^\perp$, and set
$D_\varepsilon=\varepsilon P+P^\perp$ and
$\rho_\varepsilon(\xi)=\min_{\eta\in W_k\xi}(\varepsilon|P\eta|+|P^\perp\eta|)$.
Take a ball $B$ centered at $x_*$ whose closure misses the active hyperplanes
with roots outside $J$, and $f_\varepsilon\in C_c^\infty(B)$, $0<\varepsilon\le1$,
supported in $\{|Px|\le C\varepsilon\}$ with
$\sup_\varepsilon\|f_\varepsilon\circ D_\varepsilon\|_{C^m}<\infty$ for each $m$.
Then for every integer $M>d$,
\[
 |\mathcal F_k f_\varepsilon(\xi)|
 \le C_M\varepsilon^{r+2\gamma_*}(1+\rho_\varepsilon(\xi))^{-M},
 \qquad \|\mathcal F_k f_\varepsilon\|_1\le C'_M\varepsilon^{2\gamma_*}.
\]
\end{proposition}

\begin{proof}
For $h\in C_c^\infty(B)$ supported in $K_\varepsilon=B\cap\{|Px|\le C\varepsilon\}$,
put $p_m(h)=\|h\circ D_\varepsilon\|_{C^m}$ and
$S_\xi(h)=\max_{\eta\in W_k\xi}|\mathcal F_kh(\eta)|$. The kernel bound
of \eqref{eq:FK} and the factors of $w_k$ indexed by $J$ give
\[
 S_\xi(h)\le p_0(h)\int_{K_\varepsilon}w_k\,dx
 \le C p_0(h)\,\varepsilon^{2\gamma_*}|K_\varepsilon|
 \le C_0\varepsilon^{r+2\gamma_*}p_0(h).
\]
Indeed $|\langle\alpha,x\rangle|\le C|\alpha|\varepsilon$ for $\alpha\in J$,
and $K_\varepsilon$ has normal thickness $O(\varepsilon)$ in $r$ dimensions.

For $v=D_\varepsilon e$, $|e|\le1$, let $T_v^J$ retain $\partial_v$ and the
$J$-terms of $T_v$. For the remaining roots put
$b_\alpha(x)=k_\alpha\langle\alpha,v\rangle/\langle\alpha,x\rangle$ on $B$,
extending $b_\alpha h$ by zero. Lemma~\ref{lem:dual-smooth} gives
\[
i\langle v,\eta\rangle\mathcal F_kh(\eta)=\mathcal F_k(T_vh)(\eta)=\mathcal F_k(T_v^Jh)(\eta)+\sum_{\alpha\notin J}\bigl[\mathcal F_k(b_\alpha h)(\eta)+\mathcal F_k(b_\alpha h)(s_\alpha\eta)\bigr].
\]
The reflected term has a plus sign because $1/\langle\alpha,x\rangle$ is odd
under $s_\alpha$, $w_k$ is invariant, and $E_k(s_\alpha x,s_\alpha\eta)=E_k(x,\eta)$
\cite{r18}. All tests on the right are supported in $K_\varepsilon$:
the $J$-reflections fix $x_*$ and preserve $L$. They also commute with
$D_\varepsilon$, giving $(T_v^Jh)\circ D_\varepsilon=T_e^J(h\circ D_\varepsilon)$.
Integrating the $J$-divided differences along $[s_\alpha u,u]$ and using the
positive lower bounds of the other absolute denominators on $\overline B$ gives
\[
 p_m(T_v^Jh)+\sum_{\alpha\notin J}p_m(b_\alpha h)\le C_m p_{m+1}(h)
 \qquad(|e|\le1).
\]

Let $\rho=\rho_\varepsilon(\xi)$. For $\xi\ne0$ take $\eta\in W_k\xi$
attaining $S_\xi(h)$ and $e=D_\varepsilon\eta/|D_\varepsilon\eta|$.
Then $\langle v,\eta\rangle=|D_\varepsilon\eta|\ge\rho/\sqrt2$.
For $A_m=\sup_{p_m(h)\le1}S_\xi(h)$ the multiplier and seminorm estimates
give $\rho A_{m+1}\le C_m A_m$, also valid at $\xi=0$. Since $A_{m+1}\le A_m$,
\[
 |\mathcal F_kf_\varepsilon(\xi)|
 \le p_M(f_\varepsilon)A_M
 \le C_M p_M(f_\varepsilon)(1+\rho)^{-M}A_0
 \le C'_M\varepsilon^{r+2\gamma_*}(1+\rho)^{-M}.
\]
Bounding the maximum over $W_k\xi$ by the group sum and rescaling $L$ gives
\[
 \begin{aligned}
 \|\mathcal F_k f_\varepsilon\|_1
 &\le C\varepsilon^{r+2\gamma_*}\int_\aaE(1+\rho_\varepsilon(\xi))^{-M}\,d\xi\\
 &\le C|W_k|\varepsilon^{r+2\gamma_*}
 \int_\aaE(1+\varepsilon|P\xi|+|P^\perp\xi|)^{-M}\,d\xi\\
 &= C|W_k|\varepsilon^{2\gamma_*}
    \int_{L\times L^\perp}(1+|u|+|z|)^{-M}\,du\,dz
 \ \le C''_M\varepsilon^{2\gamma_*},
 \end{aligned}
\]
since $M>d$.
\end{proof}

\begin{corollary}\label{cor:wall}
With $x_*,J,L,B$ as in Proposition~\ref{prop:dual-wall}, let
$\chi\in C_c^\infty(B)$ and $\psi\in C_c^\infty(L)$ be nonnegative with
$\psi(0)=1$, and $f_\varepsilon(x)=\chi(x)\psi(Px/\varepsilon)$. Then
$\sup_{y\in\aaE}\nu^k_y(f_\varepsilon)\le C\varepsilon^{2\gamma_*}$, and \(\nu_y^k\!\left(\bigcup_{\alpha\in R_{k,+}}\alpha^\perp\right)=0\)
for every \(y\in\aaE\).
\end{corollary}

\begin{proof}
Since $(f_\varepsilon\circ D_\varepsilon)(u)=\chi(D_\varepsilon u)\psi(Pu)$
and $f_\varepsilon|_{L^\perp}=\chi|_{L^\perp}$, positivity, Fourier inversion
and Proposition~\ref{prop:dual-wall} give
\[
 0\le\int_{L^\perp}\chi\,d\nu_y^k
 \le \nu_y^k(f_\varepsilon)
 \le (2\pi)^{-d}\|\mathcal F_k f_\varepsilon\|_1
 \le C\varepsilon^{2\gamma_*}\longrightarrow0
 \qquad(y\in\aaE).
\]
Each stratum with fixed nonempty $J$ lies in $L^\perp$ and admits a countable
cover by such patches $\{\chi>0\}$.
\end{proof}

\subsection{The nonspanning case}

This subsection carries out step~(d) of Section~\ref{sec:main}.
For $E=\spanR R_k$, write $x=x_E+x_\perp$ and $y=y_E+y_\perp$ in
$E\oplus E^\perp$. The chain fixes $x_\perp$, hence
$\mu_x^k=\mu_{x_E}^{k,E}\otimes\delta_{x_\perp}$ as in
\cite[Remark~5]{r19}, and $V_k(g_1\otimes g_2)=(V_k^Eg_1)\otimes g_2$.
Since $w_k(x)=w_{k,E}(x_E)$, \eqref{eq:rep} and ordinary Fourier inversion
give, on product tests,
\[
 \begin{aligned}
 \mathcal F_k(f_1\otimes f_2)(\xi)&=\mathcal F_{k,E}f_1(\xi_E)\,\widehat f_2(\xi_\perp),\\
 \nu_y^k(f_1\otimes f_2)&=\nu_{y_E}^{k,E}(f_1)\,f_2(y_\perp)
 =(\nu_{y_E}^{k,E}\otimes\delta_{y_\perp})(f_1\otimes f_2).
 \end{aligned}
\]
Product tests determine Radon measures. Theorem~\ref{thm:dual} on $E$ gives
\begin{equation}\label{eq:factor}
 \nu_y^k=\nu_{y_E}^{k,E}\otimes\delta_{y_\perp}\ll\lambda_{y_\perp+E},
 \qquad \nu_y^k=\delta_y\quad\text{if }E=\{0\}.
\end{equation}

\section{Moment formulas and examples}\label{sec:examples}

\subsection{Moments and power sums on a permutahedron}

Put $\iota(v)=v$, $q(v)=|v|^2$ and
$B=\sum_{\alpha\in R_+}k_\alpha(I-s_\alpha)
=2\sum_{\alpha\in R_+}k_\alpha\,(\alpha\otimes\alpha)/|\alpha|^2\ge0$.
The mean and the mean squared norm are
\[
\begin{aligned}
 (V_k\iota)(x)
 &=\int_{\aaE} y\,d\mu_x^k(y)=\E Ax
   =\int_0^\infty e^{-s}e^{-sB}x\,ds=(I+B)^{-1}x,\\
 \int_{\aaE} q\,d\mu_x^k
 &=(V_kq)(x)=\tfrac12\sum_i x_i(T_{e_i}V_kq)(x)
   =\langle x,(V_k\iota)(x)\rangle=\langle x,(I+B)^{-1}x\rangle .
\end{aligned}
\]
Indeed, \eqref{eq:chain} and Fubini give the first line since
$Q\iota=-B\iota$ implies $\E_x X_s=(e^{sQ}\iota)(x)=e^{-sB}x$.
For the second, $V_kq$ is homogeneous of degree two and $W$-invariant by
equivariance \cite{r18}, so Euler's identity applies with
$\partial_i V_kq=T_{e_i}V_kq=V_k\partial_iq=2(V_k\iota)_i$.

For $S_n$ on $\R^n$ with $n\ge2$ and constant $k\ge0$, write
$p_j(x)=\sum_i x_i^j$, $\bar x=p_1(x)/n$ and $\mathbf1=(1,\ldots,1)$.
Since $Bx=k\sum_{i<j}(I-s_{ij})x=nk(x-\bar x\mathbf1)$,
\[
 (V_kp_2)(x)=\langle x,(I+B)^{-1}x\rangle
 =\Big\langle x,\bar x\mathbf1+\frac{x-\bar x\mathbf1}{1+nk}\Big\rangle
 =\frac{p_1(x)^2}{n}+\frac{p_2(x)-p_1(x)^2/n}{1+nk}.
\]
On $K_x$, $p_1(y)=p_1(x)$, whence $V_kp_1=p_1$ and $V_k(p_1^2)=p_1^2$.
Thus $p_2-p_1^2/n$ has eigenvalue $(1+nk)^{-1}$ under $V_k$, giving
\[
 \int_{\aaE}|y-\bar x\mathbf1|^2\,d\mu_x^k(y)
 =\bigl[V_k(p_2-p_1^2/n)\bigr](x)
 =\frac{p_2(x)-p_1(x)^2/n}{1+nk}
 =\frac{|x-\bar x\mathbf1|^2}{1+nk}.
\]
For the hexagon of Figure~\ref{fig:hexagon}, $x=(2,1,0)$ and this moment is
$2/(1+3k)$.

For $k>0$, $G_s$ is irreducible on the permutation matrices, whose convex
hull is the Birkhoff polytope $B_n$. Thus $A=\int_0^\infty e^{-s}G_s\,ds$
is doubly stochastic, and Proposition~\ref{prop:simplex} and \eqref{eq:class} give
\[
 \Law(A)\sim\mathbf1_{\ri B_n}\lambda_{\aff B_n},\qquad
 \Law(A\lambda)=\mu_\lambda^k
     \sim\mathbf1_{\ri P_\lambda}\lambda_{\aff P_\lambda}.
\]
For a partition $\lambda$, $P_\lambda$ is the Newton polytope of $s_\lambda$,
so for $\lambda=(2,1,0)$ the reference measure is area on the hexagon.

\subsection{The two-point orbit}\label{sec:rank-one}

In dimension one, $s_\alpha u=-u$. For $k>0$ and $x\ne0$, set
$w(t)=(1-t)^{k-1}(1+t)^k$ on $(-1,1)$ and write $\rho_x=d\mu_x^k/dy$.
The density formula \cite[Equation~(1.2)]{r17} and the change of variables
$y=xt$ give
\[
\begin{aligned}
 \mu_x^k&=(t\mapsto xt)_*
   \bigl(C_k\mathbf1_{(-1,1)}(t)w(t)\,dt\bigr),
 \qquad C_k=\frac{\Gamma(k+\tfrac12)}{\sqrt{\pi}\,\Gamma(k)},\\
 \|\rho_x\|_1&=C_k\int_{-1}^1 w(t)\,dt=1,\\
 \|\rho_x\|_2^2&=\frac{C_k^2}{|x|}\int_{-1}^1(1-t)^{2k-2}(1+t)^{2k}\,dt
     <\infty \quad\Longleftrightarrow\quad k>\tfrac12 .
\end{aligned}
\]
The endpoint exponents are $2k>-1$ at $-1$ and $2k-2$ at $1$.
Thus the uniform multiplicity threshold in \cite[Theorem~7.2]{r11}
is sharp across root systems. For $x=0$ the law is $\delta_0$, and for
$k=0$ it is $\delta_x$.

\section*{Acknowledgments}
\addcontentsline{toc}{section}{Acknowledgments}
N.J.\ is partially supported by the Simons Foundation
(No.~\mbox{MP-TSM-00002518}). 
N.L.\ is supported by Beijing Natural Science Foundation
(No.~1264053) and the Fundamental Research Funds for the Central Universities
(No.~JZ2026HGTS0237).
\section*{Generative AI Disclosure}

The mathematical ideas and proofs in this paper are due to the authors.
ChatGPT was used to verify examples and to draft and edit the manuscript
text from the authors' notes. The authors checked every argument and
statement in the final text and take full responsibility for its content.
\begingroup
\makeatletter
\renewcommand\section{\@startsection{section}{1}{\z@}
  {.45\baselineskip}{.3\baselineskip}{\normalfont\Large\bfseries}}
\makeatother
\pdfbookmark[1]{References}{references}

\endgroup
\end{document}